\documentclass[12pt]{amsart}
\usepackage{amssymb}

\begin{document}
\newcommand{\IC}{{{\mathcal I}_C}}
\newcommand{\blub}{{{\mathcal E}nd(E)}}
\newcommand{\OP}{{\mathcal O}_{{\mathbb P}^n}}
\newcommand{\HH}{{\mathrm H}}
\newcommand{\OX}{{\mathcal O}_X}
\newcommand{\OY}{{\mathcal O}_Y}
\newcommand{\OD}{{\mathcal O}_D}
\newcommand{\ol}[1]{\overline{#1}}

\newcommand{\sta}[1]{\stackrel{#1}{\to}}
\newcommand{\ul}[1]{\underline{#1}}
\newtheorem{thm}{Theorem}
\newtheorem*{mthm}{Main~Theorem}
\newtheorem{lemma}{Lemma}
\newtheorem{cor}{Corollary}
\theoremstyle{remark}
\newtheorem{ex}{Example}
\newenvironment{rmks}{
			\trivlist \item[\hskip \labelsep{\it Remarks}:]
		      }{
		     	\endtrivlist}

\title[Arithmetically Cohen-Macaulay bundles]{Arithmetically
Cohen-Macaulay bundles on threefold hypersurfaces}
\author{N.~Mohan Kumar}
\address{Department of Mathematics, Washington University in St. Louis,
St. Louis, Missouri, 63130}
\email{kumar@wustl.edu}
\urladdr{http://www.math.wustl.edu/\~{} kumar}
\author{A.~P.~Rao}
\address{Department of Mathematics, University of Missouri-St. Louis,
St. Louis, Missouri 63121}
\email{raoa@umsl.edu}
\author{G.~V.~Ravindra}
\address{Department of Mathematics, Indian Insitute of Science,
Bangalore--560012, India}
\email{ravindra@math.iisc.ernet.in}
\date{\today}
\begin{abstract}
We prove that any rank two arithmetically Cohen-Macaulay vector bundle
on a general hypersurface of degree at least six in $\mathbb{P}^4$
must be split.
\end{abstract}
\maketitle

\section{Introduction}
This article is a continuation of the study in \cite{Helvetici} of
{\it arithmetically Cohen-Macaulay} (ACM for short) bundles on
hypersurfaces. For a detailed introduction and references, we refer
the reader to that paper. 

The main result proved here is the following
\begin{mthm}\label{main}
Fix degree $d\geq 6$. There is a non-empty Zariski open set of
hypersurfaces of degree $d$ in $\mathbb{P}^4$, none of which support
an indecomposable ACM rank two bundle.
\end{mthm}
The special case when $d=6$ was proved by Chiantini and Madonna
\cite{madonna2004}. In \cite{Helvetici}, we proved the result
analogous to the main theorem for hypersurfaces in $\mathbb{P}^5$
along with partial results on $\mathbb{P}^4$. In the current paper, we
shall summarize and use some of these partial results. We will also
use the relation between rank two ACM bundles on hypersurfaces and
Pfaffians that was observed by Beauville in \cite{Beauville} and which
we did not need in \cite{Helvetici}.

The results from \cite{Helvetici} that are important for our proof
here are paraphrased below (combining theorem 1.1 (3) and corollary
2.3 of that article).

\begin{thm}\label{previous}
Let $E$ be an indecomposable rank two ACM bundle on a smooth
hypersurface $X$ of degree $d$ in $\mathbb{P}^4$. Then
$\HH^2_*(X,\blub)$ is a non-zero cyclic module of finite length, with
the generator living in degree $-d$. If $d\geq 5$ and $X$ is general,
then $\HH^2(X,\blub) = 0$.
\end{thm}

\section{ACM bundles and Pfaffians}
We work over an algebraically closed field of characteristic zero.
Let $X\subset\mathbb{P}^n$ be a hypersurface of degree $d$. Recall
that a vector bundle $E$ of rank two on $X$ is called ACM if
$\HH^i(E(k))=0$ for all $k$ and $0<i<n-1$. By Horrocks' criterion
\cite{Horrocks}, this is equivalent to saying that $E$ has a
resolution,
\begin{equation}\label{basic1}
0\to F_1\stackrel{\Phi}\to F_0\stackrel{\sigma}\to E\to 0,
\end{equation}
where the $F_i$'s are direct sums of line bundles on $\mathbb{P}^n$.
We will assume that this resolution is minimal, with $F_0 =
\oplus_{i=1}^n \OP(-a_i)$ where $a_1\leq a_2 \leq \dots \leq
a_n$. Using \cite{Beauville}, we may write $F_1$ as $F_0^{\vee}(e-d)$,
where $e$ is the first Chern class of $E$, and we may assume that
$\Phi$ is a skew-symmetric $n\times n$ matrix with $n$ even. The
$(i,j)$-th entry $\phi_{ij}$ of $\Phi$ has degree $d-e-a_i-a_j$. The
condition of minimality implies that there are no non-zero scalar
entries in $\Phi$ and thus every degree zero entry must be zero.

We quote some facts about Pfaffians and refer the reader to
\cite{Northcott} for more details. Let $\Phi=(\phi_{ij})$ be an
$n\times n$ even-sized skew symmetric matrix and let $\mbox{Pf}(\Phi)$
denote its {\it Pfaffian}. Then $\mbox{Pf}(\Phi)^2=\det{\Phi}$. Let
$\Phi(i,j)$ be the matrix obtained from $\Phi$ by removing the $i$-th
and $j$-th rows and columns. Let $\Psi$ be the skew-symmetric matrix
of the same size with entries $\psi_{ij}=
(-1)^{i+j}\mbox{Pf}(\Phi(i,j))$ for $0\leq i<j \leq n$. We shall refer
to $\mbox{Pf}(\Phi(i,j))$ as the $(i,j)$-Pfaffian of $\Phi$.  The
product $\Phi\Psi=\mbox{Pf}(\Phi)\mbox{I}_{n}$ where $\mbox{I}_{n}$ is
the identity matrix.

\begin{ex}\label{ex}
Let $n=4$ above. Then
$$\mbox{Pf}(\Phi)=\phi_{12}\phi_{34}-\phi_{13}\phi_{24}+\phi_{14}\phi_{23}.$$
\end{ex}

The following lemma shows the relation between skew-symmetric
matrices, ACM rank $2$ bundles and the equation defining the
hypersurface.

\begin{lemma}
Let $E$ be a rank $2$ ACM bundle on a smooth hypersurface
$X\subset\mathbb{P}^4$ of degree $d$ and let $\Phi:F_1\to F_0$ be the
minimal 
skew-symmetric matrix associated to $E$. Then $X=X_{\Phi}$, the zero
locus of $\rm{Pf}(\Phi)$. Conversely, let $\Phi:F_1\to F_0$ be a
minimal skew-symmetric matrix such that the hypersurface $X_{\Phi}$
defined by $\rm{Pf}(\Phi)$ is smooth of degree $d$. Then $E_{\Phi}$,
the cokernel of 
$\Phi$, is a rank $2$ ACM bundle on $X_{\Phi}$.
\end{lemma}

\begin{proof} 
Let $f\in \HH^0(\mathbb{P}^n, \OP(d))$ be the polynomial defining
$X$. Since $E$ is supported along $X$, $\det\Phi=f^n$ for some $n$
upto a non-zero constant where $\Phi$ is as in resolution
(\ref{basic1}).  Locally $E$ is a sum of two line bundles and so the
matrix $\Phi$ is locally the diagonal matrix
$(f,f,1,\cdots,1,1)$. Since the determinant of this diagonal matrix is
$f^2$, we get $f=\mbox{Pf}(\Phi)$ (upto a non-zero constant).

To see the converse: let $\Phi$ be any skew-symmetric matrix and
$\Psi$ be  defined as above. Let $f=\mbox{Pf}(\Phi)$ be the
Pfaffian. Since $\Phi\Psi=f\mbox{I}_{n}$, this implies that the
composite $F_0(-d) \stackrel{f}{\to} F_0 \to E_{\Phi}$ is zero. Thus
$E_{\Phi}$ is 
annihilated by $f$ and so is supported on the hypersurface $X_{\Phi}$
defined by $f$. Since $X_{\Phi}$ is smooth, by the Auslander-Buchsbaum
formula, $E_{\Phi}$ is a vector bundle on $X_{\Phi}$. Therefore
locally $\Phi$ is a diagonal matrix of the form
$(f,\cdots,f,1,\cdots,1)$ where the number of $f$'s in the diagonal is
equal to the rank of $E$. Since $\det(\Phi)=f^2$, we conclude that
rank of $E_{\Phi}$ is $2$.
\end{proof}

Let $V\subset \mbox{Hom}(F_1,F_0)$ be the subspace consisting of all
minimal skew-symmetric homomorphisms, where $F_i$'s are as above. The
following is an easy consequence of the above lemma.

\begin{lemma}\label{general}
Let $\Phi_0\in V$ be an element such that $E_{\Phi_0}$ is a rank $2$
ACM bundle on a smooth hypersurface $X_{\Phi_0}$.  Then there exists a
Zariski open neighbourhood $U$ of $\Phi_0$ such that for any $\Phi\in
U$, $X_{\Phi}$ is a smooth hypersurface and $E_{\Phi}$ is a rank two
ACM bundle supported on $X_{\Phi}$.
\end{lemma}

\section{Special cases}
The proof of the main theorem will require the study of some special
cases, which are listed below.

\begin{lemma}\label{early}
Consider the following three types of curves in $\mathbb{P}^4$:
\begin{itemize}
\item a curve $C$ which is the complete intersection of three general
hypersurfaces, two of which are of degree $\leq 2$.
\item a curve $D$ which is the locus of vanishing of the principal
  $4\times 4$ 
sub-Pfaffians of a general $5\times 5$ skew-symmetric matrix $\chi$ of
linear forms. 
\item a curve $C_r$, $r\geq 0$, which is the locus of vanishing of the
$2\times 2$ minors of a general $4\times 2$ matrix $\Delta$ with one row
consisting of forms of degree $1+r$, and the remaining three rows
consisting of linear forms.
\end{itemize}
The general hypersurface $X$ in $\mathbb{P}^4$ of degree $\geq 6$
cannot contain any curve of the the first two types. The general
hypersurface $X$ of degree $d\geq\max\{6,r+4\}$ cannot contain any
curve of the third type.
\end{lemma}
 
\begin{proof}
The curve $C$ is smooth if the hypersurfaces are general. If $\chi$ is
general, the curve $D$ is smooth (see \cite{Okonek}, page 432 for
example). If $\Delta$ is general, the curve $C_r$ is smooth (see {\it
op.~cit.} page 425).

The proof of the lemma is a straightforward dimension count. By
counting the dimension of the set of all pairs $(Y,X)$ where $Y$ is a
smooth curve of the described type and $X$ is a hypersurface of degree
$d$ containing $Y$, it suffices to show that this dimension is less
than the dimension of the set of all hypersurfaces $X$ of degree $d$
in $\mathbb{P}^4$.  This can be done by showing that if $\mathcal{S}$
denotes the (irreducible) subset of the Hilbert scheme of curves in
$\mathbb{P}^4$ parameterizing all such smooth curves $Y$, then the
dimension of $\mathcal{S}$ is at most $h^0(\OY(d))-1$.

This argument was carried out in \cite{4by4} where $Y$ is any complete
intersection curve in $\mathbb{P}^4$. The case where $Y$ equals the
first type of curve $C$ in the list above is Case 2 of
\cite{4by4}. Hence we will only consider the types of curves $D$ and
$C_r$ here.

If $Y$ is of  type $D$ in the list, the sheaf $\mathcal{I}_D$ has
the following free resolution (\cite{Okonek}, page 427):
$$0 \to \OP(-5) \to \OP(-3)^{\oplus 5} \stackrel{\chi}{\to}
\OP(-2)^{\oplus 5} \to 
\mathcal{I}_D \to 0.$$  
Computing Hilbert polynomials, we see that $D$ is a smooth elliptic
quintic in $\mathbb{P}^4$, and it easily computed that that
$h^0(\mathcal N_D) = 25$. Since $h^0(\OD(d))=5d$, for
$d\geq 6$, we get $\dim{\mathcal{S}}\leq h^0(\mathcal N_D) \leq
h^0(\OD(d))-1$. 

If $Y$ is of type $C_r$ in the list, we may analyze the dimension of
the parameter space of all such $C_r$'s as follows. Let $S$ be the
cubic scroll in $\mathbb{P}^4$ given by the vanishing of the two by
two minors of the linear $3\times 2$ submatrix $\theta$ of the $4
\times 2$ matrix $\Delta$. The ideal sheaf of the determinantal
surface $S$ has resolution:
$$ 0 \to \OP(-3)^2 \stackrel{\theta}{\to} \OP(-2)^3 \to
\mathcal{I}_{S} \to 0.$$ From this one computes the dimension of the
set of such cubic scrolls to be $18$, since the $30$ dimensional space
of all $3\times 2$ linear matrices is acted on by automorphisms of
$\OP(-3)^2 $ and $\OP(-2)^3$, with scalars giving the stabilizer of
the action. Furthermore, by dualizing the resolution, we get a
resolution for $\omega_S$:
$$0 \to \OP(-5) \to \OP(-3)^{\oplus
3}\stackrel{\theta^\vee}\to\OP(-2)^{\oplus 2} \to \omega_S \to 0.$$ A
section of $\omega_S(r+3)$ gives a lift $\OP(-r-3)
\stackrel{\alpha}\to \OP(-2)^{\oplus 2}$, and we obtain a $4\times 2$
matrix
$\left(\begin{array}{c}
\theta \\
\alpha^\vee \\
\end{array}
\right)$ of the required type. Hence $C_r$ is a curve on $S$ in the
linear series $|K_S + (r+3)H|$, where $H$ is the hyperplane section on
$S$. Intersection theory on $S$ gives $K_S.K_S= 8$, $K_S.H=-5$ and
$H.H=3$. Using this, we may compute the dimension of the linear system
of $C_r$ on $S$, and we get the dimension of the set $\mathcal{S}$ of
all such $C_r$ in $\mathbb{P}^4$ to be $21$ if $r=0$ and $(3/2)r^2 +
(13/2)r +24$ otherwise.
							
The ideal sheaf of $C_r$ has a free resolution given by the
Eagon-Northcott complex \cite{EN}
\begin{multline*}
 0 \to \OP(-r-4)^{\oplus 3} \to \OP(-3)^{\oplus 2}\oplus \OP(-r-3)
^{\oplus 6} \\
\to \OP(-2)^{\oplus 3}\oplus \OP(-r-2)^{\oplus 3} \to
\mathcal{I}_{C_r} \to 0.
\end{multline*}
Let $d \geq max\{6,r+4\}$ be chosen as in the statement of the lemma.
Then $d = r+s+4$ where $s\geq 0$ ($s\geq 2$ when $r=0$; $s\geq 1$
when $r=1$). Using the above resolution, a calculation gives
$$ h^0(\mathcal O_{C_r}(d)) = \frac 32 r^2 + \frac {29}2r + 3rs + 4s
+17.$$ The required inequality $\dim{\mathcal{S}}<h^0(\mathcal
O_{C_r}(d))$ is now evident.
\end{proof}

\section{Proof of main theorem}
In this section, $E$ will be an indecomposable ACM bundle of rank two
and first Chern class $e$ on a smooth hypersurface $X$ of degree $d$
in $\mathbb{P}^4$. The minimal resolution (\ref{basic1}) gives
$\sigma: F_0 \to E \to 0$, and we may describe $\sigma$ as
$[s_1,s_2,\cdots, s_n]$ where ${s_1,s_2,\cdots, s_n}$ is a set of
minimal generators of the graded module $H_\ast^0(E)$ of global
sections of $E$, with degrees $a_1\leq a_2\leq \dots \leq a_n$.

\begin{lemma}
If $E$ is an indecomposable rank $2$ ACM bundle with first Chern class
$e$ on a general hypersurface $X\subset\mathbb{P}^4$ of degree $d\geq
6$, then there is a relation in degree $3-e$ among the minimal
generators of $S^2E$.
\end{lemma}

\begin{proof}
Consider the short exact sequence 
$$ 0 \to \OX \to \blub \to (S^2E)(-e) \to 0.$$ $S^2E(-e)$ has the same
intermediate cohomology as $\blub$ since the sequence splits in
characteristic zero.  

Choose a minimal resolution of $S^2E$: 
$$ 0 \to B \to C \to S^2E \to 0, $$ where $C$ is a direct sum of line
bundles on $X$ and $B$ is a bundle on $X$ with $\HH^1_\ast(X,B)=0$.

We first show that $B^{\vee}(e+d-5)$ is not regular. For this,
consider the dual sequence $0\to (S^2E)^\vee \to C^\vee \to B^\vee \to
0$. 

By Serre duality and Theorem \ref{previous}
$$\HH^1(X,(S^2E)^\vee(d+e-5)) = 0.$$ Therefore
$$\HH^0(X,C^\vee(d+e-5))\to\HH^0(X,B^\vee(d+e-5))$$ is onto. If
$B^\vee(d+e-5)$ were regular, the same would be true for
$$\HH^0(X,C^\vee(d+e-5+k))\to\HH^0(X,B^\vee(d+e-5+k)) \hspace{5mm}
\forall ~~k\geq 0. $$ However, this is false for $k=d$ since by Serre
duality and Theorem \ref{previous}, $\HH^1(X,(S^2E)^\vee(2d+e-5))\neq
0$. Thus $B^{\vee}(e+d-5)$ is not regular. Now
\begin{multline*}\HH^1(X,B^{\vee}(e+d-6))\cong \HH^2(X,B(1-e))\cong
\HH^1(X,S^2E(1-e))\\ \cong \HH^1(X,\blub(1)).\end{multline*} By Serre duality,
$$\HH^1(X,\blub(1))\cong \HH^2(X,\blub(d-6))$$ which by Theorem
\ref{previous} equals zero for $d\geq 6$ (this is the main place where
we use the hypothesis that $d\geq 6$). Furthermore,
$\HH^2(X,B^{\vee}(e+d-7)) = 0$ since $\HH^1_\ast(X,B)=0$. Since
$B^{\vee}(e+d-5))$ is not regular, we must have
$\HH^3(X,B^{\vee}(e+d-8)) \neq 0$.

In conclusion, $\HH^0(X,B(3-e))\neq 0$. In other words, there is a
relation in degree $3-e$ among the minimal generators of $S^2E$.
\end{proof}

\begin{lemma}\label{mu}
Let $E$ be as above. Then $1 \leq a_1+a_2+e \leq a_1+a_3+e \leq 2$.
\end{lemma}

\begin{proof}
The resolution (\ref {basic1}) for $E$ gives an exact sequence of
vector bundles on $X$: $0 \to G \to \overline F_0 \stackrel{\sigma}\to
E \to 0,$ where $\overline F_0 = F_0 \otimes \OX$ and $G$ is the
kernel.  This yields a long exact sequence,
$$0\to\wedge^2 G\to \overline{F}_0\otimes G\to S^2 \overline{F}_0\to
S^2 E\to 0.$$
From the arguments after Lemma 2.1 of \cite{Helvetici} (using
formula (5)), it follows that $\HH^2_*(\wedge^2 G) = 0$.  Hence the
map $S^2\overline F_0 \to S^2E$ is surjective on global sections. The
image of this map picks out the sections $s_is_j$ of degree $a_i+a_j$
in $S^2E$. Observe that the lowest degree minimal sections $s_1,s_2$
of $E$ induce an inclusion of sheaves $\OX(-a_1)\oplus \OX(-a_2)
\stackrel{[s_1,s_2]}\hookrightarrow E$ whose cokernel is supported on
a surface in the linear system $\vert\OX(a_1+a_2+e) \vert$ on $X$ (a
nonempty surface when $E$ is indecomposable). Hence $1 \leq a_1+a_2
+e$. There is an induced inclusion
$$S^2[\OX(-a_1)\oplus\OX(-a_2)] \hookrightarrow S^2E.$$ Therefore the
three sections of $S^2E$ given by $s_1^2, s_1s_2,s_2^2$ cannot have
any relations amongst them. Since these are also three sections of
$S^2E$ of the lowest degrees, they can be taken as part of a minimal
system of generators for $S^2E$. It follows that the relation in
degree $3-e$ among the minimal generators of $S^2E$ obtained in the
previous lemma must include minimal generators other than $s_1^2,
s_1s_2,s_2^2$. Since the other minimal generators have degree at least
$a_1+a_3$, and since we are considering a relation amongst minimal
generators, we get the inequality $a_1+a_3\leq 2-e$.
\end{proof}

\begin{lemma}\label{subpfaffians} 
For any choice of $1\leq i < j \leq n$, the $(i,j)$-Pfaffian of $\Phi$,
is non-zero.  Consequently, its degree (which is $a_i+a_j+e$) is at
least $(n-2)/2$.
\end{lemma}

\begin{proof}
On $X$, $E$ has an infinite resolution
$$\cdots \to \ol{F}_0^\vee(e-2d) \stackrel{\ol\Phi}\to \ol{F}_0(-d)
\stackrel{\ol\Psi}\to \ol{F}_0^\vee(e-d) \stackrel{\ol\Phi}\to \ol{F}_0 \to
E \to 0.$$
We also have
\[
\begin{array}{cccccccccc}
& \ol{F}_0^\vee(e-2d) & \stackrel{\ol\Phi}\to &
\ol{F}_0(-d) & \stackrel{\sigma}\to & E(-d) & & & &  \\ 

    &                     &                        &
                  &                       & 
\downarrow{\stackrel{\alpha}\cong}               & & & &  \\
                                                   
    &                     &                        & 
                  &                       &  E^\vee(e-d) & 
\stackrel{\sigma^\vee}\hookrightarrow &
\ol{F}_0^\vee(e-d) & \stackrel{\ol{\Phi}^\vee}\to & \ol{F}_0. \\
\end{array}
\]
Let $\ol{\Theta}=\sigma^\vee\alpha\sigma$. Since
$\sigma=(s_1,\cdots,s_n)$, we may express the $(i,j)$-th entry of
$\ol{\Theta}$ as $\theta_{ij}=s_i^\vee s_j$ (suppressing the canonical
isomorphism $\alpha$). $\Phi^\vee=-\Phi$ and
$\alpha\sigma:\ol{F}_0(-d) \to E^\vee(e-d)$ is surjective on global
sections. Hence we have a commuting diagram 
\[
\begin{array}{cccccccccccc}
\ol{F}_0^\vee(e-2d) & \stackrel{\ol\Phi}\to & \ol{F}_0(-d) &
\stackrel{\ol\Psi}\to & \ol{F}_0^\vee(e-d) & \stackrel{\ol\Phi}\to &
\ol{F}_0 & \to &  E & \to & 0 \\
& & \downarrow{\stackrel{B}\cong} & &
\downarrow{\stackrel{-I}\cong} & & || & & || & & \\
\ol{F}_0^\vee(e-2d) & \stackrel{\ol\Phi}\to & \ol{F}_0(-d)&
\stackrel{\ol\Theta}\to & \ol{F}_0^\vee(e-d) & \stackrel{-\ol\Phi}\to & 
\ol{F}_0 & \to & E & \to & 0 \\
\end{array}
\]
It is easy to see that $B$ is an isomorphism. As a result, every
column of $B$ has a non-zero scalar entry. 

Now suppose that $\ol\psi_{ij}=0$ for some $i,~j$ so that $\sum_k
s_i^\vee s_k b_{kj}=0$. Let $Y_i$ be the curve given by the vanishing
of the minimal section $s_i$ with the exact sequence
$$0 \to \OX(-a_i) \stackrel{s_i}\to E \stackrel{s_i^\vee}\to
I_{Y_i/X}(a_i+e) \to 0.$$
Hence $s_i^\vee s_i=0$ and $s_i^\vee s_k$ for $k\neq i$ give minimal
generators for $I_{Y_i/X}$. It follows that no $b_{kj}$ can be a
non-zero scalar for $k\neq i$. Hence $b_{ij}$ has to be a non-zero
scalar and the only one in the $j$-th column. However,
$\ol\psi_{jj}=0$. So by the same argument, $b_{jj}$ is the only
non-zero scalar. To avoid contradiction, $\ol\psi_{ij}$ and hence
$\psi_{ij}\neq 0$ for $i\neq j$.
\end{proof}

We now complete the proof of the Main Theorem. As in the previous
lemmas, assume that $X$ is general of degree $d\geq 6$, with $E$ an
indecomposable rank two ACM bundle on $X$. We will show that the
inequalities of Lemma \ref{mu} lead us to the special cases of Lemma
\ref{early}, giving a contradiction.

Let $\mu = a_1+a_2+e$. By Lemma \ref{mu}, $1 \leq \mu \leq 2$.

\subsection*{Case $\mu = 1$} 
In this case, in order for the $(1,2)$-Pfaffian of $\Phi$ to be
linear, by Lemma \ref{subpfaffians}, $n$ must equal 4.  In the $4
\times 4$ matrix $\Phi$, the $(1,2)$-Pfaffian is the entry $\phi_{3
4}$ which we are claiming is linear. Likewise the $(1,3)$-Pfaffian is
the entry $\phi_{2 4}$ which by Lemma \ref{mu} has degree $a_1+a_3+e
\leq 2$.  By Lemma \ref{general}, we may assume that $\phi_{1
4},\phi_{2 4},\phi_{3 4}$ define a smooth complete intersection curve
and $X$ contains this curve by example \ref{ex}.  By Lemma
\ref{early}, $X$ cannot be general.

\subsection*{Case $\mu = 2$}
In this case $a_2=a_3$. By Lemma \ref{subpfaffians}, $n$
must be 4 or 6.  The case $n=4$ is ruled out again by the arguments of
the above paragraph since $\Phi$ has two entries of degree $2$ in its
last column.  We will therefore assume that $n=6$.
The matrix \[\Phi =
\left(
\begin{array}{cccccc}
0 & \phi_{12} & \phi_{13}&\phi_{14} & \phi_{15} & \phi_{16} \\
\ast & 0 &\phi_{23} & \phi_{24} &\phi_{25} & \phi_{26}  \\
\ast &\ast  & 0  &\phi_{34} & \phi_{35} & \phi_{36}  \\
\ast &\ast  &\ast & 0 & \phi_{45} & \phi_{46}  \\
\ast & \ast &\ast &\ast & 0 & \phi_{56} \\
\ast &\ast  & \ast&\ast  & \ast & 0  \\

\end{array}
\right)
\]
is skew-symmetric and by our choice of ordering of
the $a_i$'s, the degrees of the upper triangular entries
are non-increasing as we move to the right or down.

As remarked before, the degree of $\phi_{ij}$ is $d-e-a_i-a_j$. The
$(1,2)$-Pfaffian (which is a non-zero quadric when $\mu =2$) is given
by the expression (see example \ref{ex})
\begin{equation}\label{12pfaff}
\mbox{Pf}(\Phi(1,2))=\phi_{34}\phi_{56}-\phi_{35}\phi_{46}+\phi_{36}\phi_{45}.
\end{equation}
We shall consider the following two sub-cases, one where $\phi_{56}$
has positive degree (and hence can be chosen non-zero by Lemma
\ref{general}) and the other where it has non-positive degree (and
hence is forced to be zero):

\subsection*{}$\ul{d-e-a_5-a_6 >0}$.
\vspace{2mm}

\noindent Since $\phi_{34}\cdot \phi_{56}$ is one term in the
$(1,2)$-Pfaffian of $\Phi$, and since degree $\phi_{34}$ is at least
degree $\phi_{56}$, they are both forced to be linear. Therefore
$\phi_{34}$, $\phi_{35}$, $\phi_{36}$ have the same degree (=$1$) and
so $a_4=a_5=a_6$. Likewise, $a_3=a_4=a_5$. Therefore
$a_2=a_3=a_4=a_5=a_6$. Hence $\Phi$ has a principal $5\times 5$
submatrix $\chi$ (obtained by deleting the first row and column in
$\Phi$) which is a skew symmetric matrix of linear terms,
while its first row and first column have entries of degree $1+r,
r\geq 0$.

By Lemma \ref{general} we may assume that the ideal of the $4\times 4$
Pfaffians of $\chi$ defines a smooth curve $C$. $X$ is then a degree
$d = 3 +r $ hypersurface containing $C$. By Lemma \ref{early}, $X$
cannot be general when $d\geq 6$.

\subsection*{}$\ul{d-e-a_5-a_6 \leq 0}$.
\vspace{2mm}

\noindent In this case, the entry $\phi_{56} =0$.  Suppose $\phi_{46}$
is also zero. Then both $\phi_{36}$ and $\phi_{45}$ must be linear and
non-zero since the $(1,2)$-Pfaffian of $\Phi$ (see equation
\ref{12pfaff}) is a non-zero quadric.  Since $a_2=a_3$, $\phi_{26}$ is
also linear.  Thus using Lemma \ref{general}, $X$ contains the complete
intersection curve given by the vanishing of $\phi_{16}$ and the two
linear forms $\phi_{36},\phi_{26}$. By Lemma \ref{early}, $X$ cannot
be general.

So we may assume that $\phi_{46} \neq 0$. Since $\phi_{35}$ is also
non-zero, both must be linear. Hence $a_3+a_5=a_4+a_6$, and so
$a_3=a_4$ and $a_5=a_6$.

After twisting $E$ by a line bundle, we may assume that
$a_2=a_3=a_4=0 \leq a_5=a_6=b$. The linearity of the entry
$\phi_{46}$ gives $d-e-b=1$.  The condition $d-e-a_5-a_6 \leq 0$
yields $1\leq b$. Taking first Chern classes in resolution
(\ref{basic1}) gives $e = 2-a_1$.

Let $r = -a_1,~ s= b-1$. Then $r,~s\geq0$, and $d= r +s +4$. If we
inspect the matrix $\Phi$, the non-zero rows in columns 5 and 6 give a
$4\times 2$ matrix $\Delta$ with top row of degree $1+r$ and the other
entries all linear. By Lemma \ref{general}, we may assume that the
$2\times 2$ minors of this $4\times 2$ matrix define a smooth curve
$C_r$ as described in Lemma \ref{early}.  Since $X$ contains this
curve, $X$ cannot be general when $d\geq 6$.

\end{document}